\documentclass[12pt]{article}
\usepackage{amssymb}
\usepackage{amsfonts}

\newtheorem{pro}{Proposition}

\usepackage[bookmarksopen=true,colorlinks=true,linkcolor=blue,citecolor=blue, urlcolor=blue, pdfstartview={FitBH}, pdfview={FitBH}]{hyperref}
\begin{document}
\author{Mustapha CHELLALI  \& Alain SALINIER}
\title{La conjecture de Casas Alvero pour les degr\'{e} $5p^{e}$}
\date{}
\maketitle
\vspace*{5mm}
\begin{flushleft}
\begin{tabular}{ll}
\cline{2-2}
\ \\

{\bf \scriptsize  R\'{e}sum\'{e} :}&\scriptsize  {Selon la conjecture de Casas Alvero,  si un polyn\^{o}me \`{a} une variable  de degr\'{e} $n$ sur un corps commutatif de } \vspace{0.2cm}\\
&\begin{minipage}{15cm}
  {\scriptsize caract\'{e}ristique 0  est non premier avec chacune de ses $n-1$ premi\`{e}res d\'{e}riv\'{e}s, alors il est de forme $c(X-r)^{n}$. Soient $p$ un nombre premier et $e$ un entier,  la conjecture a \'{e}t\'{e} d\'{e}montr\'{e}e pour les polyn\^{o}mes de degr\'{e} $p^{e},2p^{e}, 3p^{e}\\ (p\neq 2) \ et \ 4p^{e}\ (p\neq 3,5,7)$. Dans ce travail on montre que la conjecture est vrai   pour les polyn\^{o}mes de degr\'{e} $5p^{e}\ (p\neq 2,3,7,11,131,193,599,3541,8009)$. On corrige aussi une erreur dans \cite{Dra-Jong} pour les degr\'{e} $4p^{e}$ }
\end{minipage} \\
\  \\
\cline{2-2}
\end{tabular}
\end{flushleft}

\section{Introduction} 

Soit $k$ un corps commutatif, $P\in k[X]$ un polyn\^{o}me de degr\'{e} $n$, la conjecture de Casas-Alvero veut que si $P$ est non premier avec chacune de ses $n-1$ premi\`{e}res d\'{e}riv\'{e}es, alors  c'est un mon\^{o}me, c'est-\'{a}-dire de la forme $c(X-r)^{n}$. Cette conjecture est \'{e}videment fausse si la caract\'{e}ristique de $k$ est $\neq 0$ comme le montre l'exemple $P=X^{p+1}-X^{p}$ en caract\'{e}ristique $p$. Suivant \cite{Hans}, il importe de modifier l'\'{e}nonc\'{e} de la conjecture en caract\'{e}ristique $\neq 0$, en rempla\c{c}ant les d\'{e}riv\'{e}s ordinaire $P^{(i)}$ par les d\'{e}riv\'{e}s de Hasse d\'{e}finies par $P(X+h)=\sum_{i}^{}P_{i}(h)X^{i}$ (en caract\'{e}ristique 0 on a simplement $P_{i}=P^{(i)}/i!$). Un lien int\'{e}ressant a \'{e}t\'{e}  trouv\'{e} entre la conjecture en caract\'{e}ristique $0$ et la conjecture en caract\'{e}ristique $p\neq 0$ (cf  \cite{Hans}  voir aussi \cite{Dra-Jong} pour une preuve \'{e}l\'{e}mentaire), en notons   $\overline{ \mathbb{F}}_{p}$ la cl\^{o}ture alg\'{e}brique de $\mathbb{F}_{p}$, il s'\'{e}nonce comme suit~:

\begin{pro}\label{} Soit $n$ un entier $\geq 1$
\begin{itemize}
\item Si pour un nombre premier donn\'{e} $p$ la conjecture de Casas-Alvero est vrai sur $\overline{ \mathbb{F}}_{p}$  pour tout polyn\^{o}me de degr\'{e} $n$, alors elle est vrai en caract\'{e}ristique $0$ pour tout polyn\^{o}me de degr\'{e}  de la forme $np^{e}\ (e\in \mathbb{N}) $ 
\item  Inversement si la conjecture est vrai en caract\'{e}ristique $0$ pour tout polyn\^{o}me de degr\'{e}   $n$, alors elle est vrai sur $\overline{ \mathbb{F}}_{p}$   pour tout polyn\^{o}me de degr\'{e} $n$, sauf peut \^{e}tre pour un nombre fini de $p$.

\end{itemize}

\end{pro}

 Comme cons\'{e}quence on obtient~:

\begin{itemize}
\item Pour $n=1$, il n'y a aucun mauvais nombre premier, par cons\'{e}quent   la conjecture est vrai pour tout polyn\^{o}me de degr\'{e} $p^{e}\ (p \ premier)$
\item Pour $n=2$,  si un polyn\^{o}me  de degr\'{e} 2 v\'{e}rifie   les hypoth\`{e}ses de Casas-Alvero, comme il a une racine double, c'est un mon\^{o}me, par suite la conjecture est vrai pour tout polyn\^{o}me de degr\'{e} $2p^{e}\ (p\  premier)$
\item  Pour $n=3$, cherchons les mauvais nombres premiers,  soit $p$ un nombre premier , soit $P$ un polyn\^{o}me de degr\'{e} 3 sur $\overline{ \mathbb{F}}_{p}$ v\'{e}rifiant les hypoth\`{e}ses de Casas-Alvero , apr\`{e}s transformation affine on se ram\`{e}ne \`{a} $P$ de la forme $X^{3}-aX^{2}$, \'{e}cartons d'abord les cas $p=2,3$ puisque la conjecture de Casas-Alvero est v\'{e}rifi\'{e}e pour $3^{e}$ et que on sait que $X^{3}-X^{2}$ est un contre exemple en caract\'{e}ristique 2

\[\left \{ \begin{array}{l}
P'=3X^{2}-2aX\\
P_{2}=3X-a

\end{array}\right.\]  

Si 0 est racine de $P''\longrightarrow a=0\longrightarrow P=X^{3}$, sinon la racine de $P''$ en commun avec $P$ est $a\longrightarrow 2a=0\longrightarrow a=0$. Ainsi $2$ est  le seul  mauvais nombre premier. 

Conclusion~: Si $p$ est nombre premier $\neq 2$ la conjecture de Casas-Alvero est vrai pour les polyn\^{o}mes de degr\'{e} $3p^{e}$

\item  Pour $n=4$, cherchons les mauvais nombres premiers,  soit $p$ un nombre premier , soit $P$ un polyn\^{o}me de degr\'{e} 4 sur $\overline{ \mathbb{F}}_{p}$ v\'{e}rifiant les hypoth\`{e}ses de Casas-Alvero , apr\`{e}s transformation affine on se ram\`{e}ne \`{a} $P$ de la forme $X^{4}-aX^{3}+bX^{2}$, \'{e}cartons d'abord les cas $p=2,3$ puisque la conjecture de Casas-Alvero est v\'{e}rifi\'{e}e pour $2^{e}$ et que on sait que $X^{4}-X^{3}$ est un contre exemple en caract\'{e}ristique 3.

\[\left \{ \begin{array}{l}
P'=4X^{3}-3aX^{2}+2bX\\
P_{2}=6X^{2}-3aX+b\\
P_{3}=4X-a  

\end{array}\right.\]  

Supposons d'abord que 0 n'est pas racine de $P''$ et $P'''$ soit $ab\neq 0$, \'{e}crivons $P=X^{2}(X-x_{1})(X-x_{2})$, soient $\alpha ,\beta$ les racines de  $P''$ et $P'''$ en commun avec $P$, on peut supposer que $\beta = x_{1}$, deux cas sont alors possibles $(\alpha ,\beta)=(x_{1},x_{1})$  ou  $(\alpha ,\beta)=(x_{2},x_{1})$
\begin{itemize}
\item  1 er cas : $(\alpha ,\beta)=(x_{1},x_{1})$

\[\left \{ \begin{array}{l}
4x_{1}=a=x_{1}+x_{2}\longrightarrow x_{2}=3x_{1}\\
6x_{1}^{2}-12x_{1}^{2}+b=0\longrightarrow b=6x_{1}^{2}\ soit\   6x_{1}^{2}=x_{1}x_{2}=3x_{1}^{2}\longrightarrow x_{1}=x_{2}=0

\end{array}\right.\]
\item 2 \`{e}me cas :  $(\alpha ,\beta)=(x_{2},x_{1})$

\[\left \{ \begin{array}{l}
4x_{1}=a=x_{1}+x_{2}\longrightarrow x_{2}=3x_{1}\\
6x_{2}^{2}-12x_{1}x_{2}+b=0\longrightarrow b=-2x_{2}^{2}\ soit\   -2x_{2}^{2}=x_{1}x_{2}=3x_{1}^{2}\longrightarrow 21x_{1}^{2}=0\longrightarrow p=7

\end{array}\right.\]

Inversement si  $p=7$, prenons $x_{1}=1$ et $x_{2}=3x_{1}=3$ et $P=X^{2}(X-1)(X-3)$ on a $P'''(1)=0$ et $P''(3)=0$, donc $P$ v\'{e}rifie les hypoth\`{e}ses de Casas-Alvero dans $\overline{ \mathbb{F}_{7}}$ et n'est pas mon\^{o}me, ainsi $7$ est un mauvais nombre premier. 
\end{itemize}

\end{itemize}  

Supposons maintenant  que 0 est  racine de $P''$ ou $P'''$ soit $a\ ou\ b=0$ (mais pas les deux car dans ce cas $P$ serait mon\^{o}me et $p$ un bon nombre premier )
\begin{itemize}
\item $a=0\longrightarrow x_{1}+x_{2}=0$, on peut supposer que la racine de $P''$ est alors $x_{1}\neq 0\longrightarrow 6x_{1}^{2}+b=0\longrightarrow b=-6x_{1}^{2}=x_{1}x_{2}=-x_{1}^{2}\longrightarrow 5x_{1}^{2}=0\longrightarrow p=5$. Inversement si $p=5$, posons $x_{1}=1$ et $x_{2}=-1$ et $P=X^{2}(X-1)(X+1)$ on a $P'(0)=P'''(0)=0$ et $P''(1)=6-1=0$, donc $P$ v\'{e}rifie les hypoth\`{e}ses de Casas-Alvero dans $\overline{ \mathbb{F}}_{5}$ et n'est pas mon\^{o}me, ainsi $5$ est un mauvais nombre premier.    Ce cas est pass\'{e} inaper\c{c}u dans \cite{Dra-Jong} suite a une erreur de consid\'{e}ration de d\'{e}terminant (cf page 33)

\item  $b=0\longrightarrow x_{1}x_{2}=0$, on peut supposer que la racine de $P'''$ est alors $x_{1}\neq 0\longrightarrow 4x_{1}-a=0\longrightarrow a=4x_{1}=x_{1}+x_{2}\longrightarrow 3x_{1}^{2}=0\longrightarrow x_{1}=x_{2}=0$. 
\end{itemize}

Conclusion~: Si $p$ est nombre premier $\neq 3,5,7$ la conjecture de Casas-Alvero est vrai pour les polyn\^{o}mes de degr\'{e} $4p^{e}$

\section{Cas des degr\'{e}s $5p^{e}$ }

Notre r\'{e}sultat principal est~:

\begin{pro}\label{} La conjecture de Casas-Alvero est vrai pour les polyn\^{o}mes de degr\'{e} $5p^{e}$, \\$e$ entier et  $p$ premier $\neq 2,3,7,11,131,193,599,3541,8009$

\end{pro}

Preuve~: Nous allons poursuivre les m\'{e}thodes ci-dessus pour d\'{e}terminer les mauvais nombres premiers pour les polyn\^{o}mes de degr\'{e} $5$. Soit $P$ un polyn\^{o}me de degr\'{e} 5 sur $\overline{ \mathbb{F}}_{p}$ v\'{e}rifiant les hypoth\`{e}ses de Casas-Alvero , apr\`{e}s transformation affine on se ram\`{e}ne \`{a} $P$ de la forme $X^{5}-aX^{4}+bX^{3}-cX^{2}$, \'{e}cartons d'abord les cas $p=2,3,5$ puisque la conjecture de Casas-Alvero est v\'{e}rifi\'{e}e pour $5^{e}$ et que $X^{5}-X^{4}$ est un contre exemple en caract\'{e}ristique 2 et $X^{5}+X^{4}$ est un contre exemple en caract\'{e}ristique 3.

\[\left \{ \begin{array}{l}
P'=5X^{4}-4aX^{3}+3bX^{2}-2cX\\
P_{2}=10X^{3}-6aX^{2}+3bX-c\\
P_{3}=10X^{2}-4aX+b\\  
P_{4}=5X-a

\end{array}\right.\]  

Soient $\alpha, \beta, \gamma  $ les racines de $P'',P''',P^{(4)}$ en commun avec $P$, nous distinguons deux cas~:

\begin{itemize}
\item  1 er cas~:  $0\notin \left \{\alpha, \beta, \gamma \right \}$ 

Ecrivons $P=X^{2}(X-x_{1})(X-x_{2})(X-x_{3})$, on peut supposer $\gamma =x_{1}$ 

\begin{itemize}
\item Cas $(\alpha, \beta, \gamma) = (x_{1},x_{1},x_{1})$

\[\left \{ \begin{array}{l}
P_{4}(x_{1})=5x_{1}-a=0\longrightarrow 5x_{1}=a=x_{1}+x_{2}+x_{3}\longrightarrow 4x_{1}=x_{2}+x_{3}\\  
P_{3}(x_{1})=0\longrightarrow 10x_{1}^{2}-20x_{1}^{2}+b=0\\\longrightarrow b=10x_{1}^{2}=x_{1}x_{2}+x_{1}x_{3}+x_{2}x_{3}=4x_{1}^{2}+x_{2}(4x_{1}-x_{2})\longrightarrow \fbox{$6x_{1}^{2}-4x_{1}x_{2}+x_{2}^{2}=0$}\\
P_{2}(x_{1})=0\longrightarrow 10x_{1}^{3}-30x_{1}^{3}+30x_{1}^{3}-x_{1}x_{2}x_{3}=0\longrightarrow 10x_{1}^{2}=x_{2}x_{3}=x_{2}(4x_{1}-x_{2})\\
\longrightarrow  \fbox{$10x_{1}^{2}-4x_{1}x_{2}+x_{2}^{2}=0$}
\end{array}\right.\]

Comme $x_{1}$ est suppos\'{e} $\neq 0$, posons $r=x_{2}/x_{1}$, on a le syst\`{e}me~:

\[\left \{ \begin{array}{l}
r^{2}-4r+6=0\\  
r^{2}-4r+10=0
\end{array}\right.\]

Le resultant de ces deux \'{e}quations est $16\neq 0$ puisque $p\neq 2$, par suite ce cas est impossible.

\item  Cas $(\alpha, \beta, \gamma) = (x_{3},x_{2},x_{1})$

\[\left \{ \begin{array}{l}
P_{4}(x_{1})=5x_{1}-a=0\longrightarrow 5x_{1}=a=x_{1}+x_{2}+x_{3}\longrightarrow 4x_{1}=x_{2}+x_{3}\\  
P_{3}(x_{2})=0\longrightarrow 10x_{2}^{2}-20x_{1}x_{2}+b=0\\\longrightarrow b=-10x_{2}^{2}+20x_{1}x_{2}=x_{1}x_{2}+x_{1}x_{3}+x_{2}x_{3}=4x_{1}^{2}+x_{2}(4x_{1}-x_{2})\\\longrightarrow \fbox{$4x_{1}^{2}-16x_{1}x_{2}+9x_{2}^{2}=0$}\\
P_{2}(x_{3})=0\longrightarrow 10x_{3}^{3}-30x_{1}x_{3}^{2}+3(-10x_{2}^{2}+20x_{1}x_{2})x_{3}-x_{1}x_{2}x_{3}=0\\\longrightarrow 10x_{3}^{2}-30x_{1}x_{3}+3(-10x_{2}^{2}+20x_{1}x_{2})-x_{1}x_{2}=0\\\longrightarrow 10(4x_{1}-x_{2})^{2}-30x_{1}(4x_{1}-x_{2})+3(-10x_{2}^{2}+20x_{1}x_{2})-x_{1}x_{2}=0\\\
\longrightarrow  \fbox{$40x_{1}^{2}+9x_{1}x_{2}-20x_{2}^{2}=0$}
\end{array}\right.\]

Comme $x_{1}$ est suppos\'{e} $\neq 0$, posons $r=x_{2}/x_{1}$, on a le syst\`{e}me~:

\[\left \{ \begin{array}{l}
9r^{2}-16r+4=0\\  
-20r^{2}+9r+40=0
\end{array}\right.\]

Le resultant de ces deux \'{e}quations est $32036=2^{2}.8009$ puisque $p\neq 2$ cela n'est possible que si $p=8009$. Inversement si  $p=8009$ on va remonter ces \'{e}quations pour construir un contre exemple modulo $p$. eliminon $r$ entre ces deux equations on obtient $r =  440/239 = 2113 \ mod \ 8009$. En  fixant $x_{1}=1$ alors $x_{2}=r$ et $x_{3}=4x_{1}-x_{2}=4-r$, cela doone le contre exemple~:

\[P=x^{2}(x-1)(x-r)(x-4+r)=x ^{5} - 5 x^{4}  - 3309 x^{3}  + 3313 x^{2}\qquad mod      \ 8009\]

On v\'{e}rifie bien que modulo 8009 on a

\[\begin{array}{r}
P_{4}(1)=0\\
P_{3}(r)=0\\
P_{2}(4-r)=0
\end{array}\]

 \item Cas  $(\alpha, \beta, \gamma) = (x_{2},x_{1},x_{1})$

\[\left \{ \begin{array}{l}
P_{4}(x_{1})=5x_{1}-a=0\longrightarrow 5x_{1}=a=x_{1}+x_{2}+x_{3}\longrightarrow 4x_{1}=x_{2}+x_{3}\\  
P_{3}(x_{1})=0\longrightarrow 10x_{1}^{2}-20x_{1}^{2}+b=0\\\longrightarrow b=10x_{1}^{2}=x_{1}x_{2}+x_{1}x_{3}+x_{2}x_{3}=4x_{1}^{2}+x_{2}(4x_{1}-x_{2})\longrightarrow \fbox{$6x_{1}^{2}-4x_{1}x_{2}+x_{2}^{2}=0$}\\
P_{2}(x_{2})=0\longrightarrow 10x_{2}^{3}-30x_{1}x_{2}^{2}+30x_{1}^{2}x_{2}-x_{1}x_{2}x_{3}=0\\\longrightarrow 10x_{2}^{2}-30x_{1}x_{2}+30x_{1}^{2}-x_{1}x_{3}=0\\\longrightarrow 10x_{2}^{2}-30x_{1}x_{2}+30x_{1}^{2}-x_{1}(4x_{1}-x_{2})=0\\\
\longrightarrow  \fbox{$26x_{1}^{2}- 29x_{1}x_{2}+10x_{2}^{2}=0$}
\end{array}\right.\]

Comme $x_{1}$ est suppos\'{e} $\neq 0$, posons $r=x_{2}/x_{1}$, on a le syst\`{e}me~:

\[\left \{ \begin{array}{l}
r^{2}-4r+6=0\\  
10r^{2}-29r+26=0
\end{array}\right.\]

Le resultant de ces deux \'{e}quations est $386=2.193$ puisque $p\neq 2$ cela n'est possible que si $p=193$. Inversement si  $p=193$ on va remonter ces \'{e}quations pour construir un contre exemple modulo $p$. eliminon $r$ entre ces deux equations on obtient $r =  34/11 = 161 \ mod \ 193$. En  fixant $x_{1}=1$ alors $x_{2}=r$ et $x_{3}=4x_{1}-x_{2}=4-r$, cela doone le contre exemple~:

\[P=x^{2}(x-1)(x-r)(x-4+r)=x^{5}  - 5x^{4}  + 10x^{3}  - 6x^{2}\qquad mod      \ 193\]

On v\'{e}rifie bien que modulo 193 on a

\[\begin{array}{r}
P_{4}(1)=0\\
P_{3}(r)=0\\
P_{2}(4-r)=0
\end{array}\]

\item  Cas $(\alpha, \beta, \gamma) = (x_{1},x_{2},x_{1})$

\[\left \{ \begin{array}{l}
P_{4}(x_{1})=5x_{1}-a=0\longrightarrow 5x_{1}=a=x_{1}+x_{2}+x_{3}\longrightarrow 4x_{1}=x_{2}+x_{3}\\  
P_{3}(x_{2})=0\longrightarrow 10x_{2}^{2}-20x_{1}x_{2}+b=0\\\longrightarrow b=-10x_{2}^{2}+20x_{1}x_{2}=x_{1}x_{2}+x_{1}x_{3}+x_{2}x_{3}=4x_{1}^{2}+x_{2}(4x_{1}-x_{2})\\\longrightarrow \fbox{$4x_{1}^{2}-16x_{1}x_{2}+9x_{2}^{2}=0$}\\
P_{2}(x_{1})=0\longrightarrow 10x_{1}^{3}-30x_{1}x_{1}^{2}+3(-10x_{2}^{2}+20x_{1}x_{2})x_{1}-x_{1}x_{2}x_{3}=0\\\longrightarrow 10x_{1}^{2}-30x_{1}x_{1}+3(-10x_{2}^{2}+20x_{1}x_{2})-x_{2}x_{3}=0\\\longrightarrow 10x_{1}^{2}-30x_{1}x_{1}+3(-10x_{2}^{2}+20x_{1}x_{2})-x_{2}(4x_{1}-x_{2})=0\\\
\longrightarrow  \fbox{$-20x_{1}^{2}+56x_{1}x_{2}-29x_{2}^{2}=0$}
\end{array}\right.\]

Comme $x_{1}$ est suppos\'{e} $\neq 0$, posons $r=x_{2}/x_{1}$, on a le syst\`{e}me~:

\[\left \{ \begin{array}{l}
9r^{2}-16r+4=0\\  
29r^{2}-56r+20=0
\end{array}\right.\]

Le resultant de ces deux \'{e}quations est $256=2^{8}$ puisque $p\neq 2$ ce cas est impossible.

\item  Cas $(\alpha, \beta, \gamma) = (x_{2},x_{2},x_{1})$

\[\left \{ \begin{array}{l}
P_{4}(x_{1})=5x_{1}-a=0\longrightarrow 5x_{1}=a=x_{1}+x_{2}+x_{3}\longrightarrow 4x_{1}=x_{2}+x_{3}\\  
P_{3}(x_{2})=0\longrightarrow 10x_{2}^{2}-20x_{1}x_{2}+b=0\\\longrightarrow b=-10x_{2}^{2}+20x_{1}x_{2}=x_{1}x_{2}+x_{1}x_{3}+x_{2}x_{3}=4x_{1}^{2}+x_{2}(4x_{1}-x_{2})\\\longrightarrow \fbox{$4x_{1}^{2}-16x_{1}x_{2}+9x_{2}^{2}=0$}\\
P_{2}(x_{2})=0\longrightarrow 10x_{2}^{3}-30x_{1}x_{2}^{2}+3(-10x_{2}^{2}+20x_{1}x_{2})x_{2}-x_{1}x_{2}x_{3}=0\\\longrightarrow 10x_{2}^{2}-30x_{1}x_{2}+3(-10x_{2}^{2}+20x_{1}x_{2})-x_{1}x_{3}=0\\\longrightarrow 10x_{2}^{2}-30x_{1}x_{2}+3(-10x_{2}^{2}+20x_{1}x_{2})-x_{1}(4x_{1}-x_{2})=0\\\
\longrightarrow  \fbox{$-4x_{1}^{2}+31x_{1}x_{2}-20x_{2}^{2}=0$}
\end{array}\right.\]

Comme $x_{1}$ est suppos\'{e} $\neq 0$, posons $r=x_{2}/x_{1}$, on a le syst\`{e}me~:

\[\left \{ \begin{array}{l}
9r^{2}-16r+4=0\\  
20r^{2}-31r+4=0
\end{array}\right.\]

Le resultant  de ces deux \'{e}quations est $-524=-2^{2}.131$ puisque $p\neq 2$ cela n'est possible que si $p=131$. Inversement si  $p=131$ on va remonter ces \'{e}quations pour construir un contre exemple modulo $p$. eliminon $r$ entre ces deux equations on obtient $r =  44/41 = 49 \ mod \ 131$. En  fixant $x_{1}=1$ alors $x_{2}=r$ et $x_{3}=4x_{1}-x_{2}=4-r$, cela doone le contre exemple~:

\[P=x^{2}(x-1)(x-r)(x-4+r)=x^{5}  - 5 x^{4}  + 26 x^{3}  - 22 x^{2}\qquad mod      \ 131\]

On v\'{e}rifie bien que modulo 131 on a

\[\begin{array}{r}
P_{4}(1)=0\\
P_{3}(r)=0\\
P_{2}(4-r)=0
\end{array}\]

\end{itemize}

\item  2 \`{e}me cas~:  $0\in  \left \{\alpha, \beta, \gamma \right \}$ 

Autrement dit $a\ ou\ b\ ou\ c=0$, on a les  resultants

\[Res_{x}(P,P_{2})=100\, a^3\, c^4 - 24\, a^2\, b^2\, c^3 - 459\, a\, b\, c^4 + 98\
 \, b^3\, c^3 + 729\, c^5\]

\[Res_{x}(P,P_{3})=96\, a^3\, b^2\, c - 18\, a^2\, b^4 - 480\, a\, b^3\, c + 81\, b\
^5 + 1000\, b^2\, c^2\]

\[Res_{x}(P,P^{(4)})=4\, a^5 - 25\, b\, a^3 + 125\, c\, a^2\]

\begin{itemize}
\item 1er cas ~:   $a=0$ 

Dans ce cas les resultants  ci-dessus deviennent~:

\[\begin{array}{l}
Res_{x}(P,P_{2})=c^3(98\, b^3\,  + 729\, c^2)\\
Res_{x}(P,P_{3})=b^2(81\, b^3 + 1000\, \, c^2)

\end{array}\]

Si $P$ v\'{e}rifie les hypoth\`{e}ses de Casas-Alvero on aura alors

\[\begin{array}{l}
c^3(98\, b^3\,  + 729\, c^2)=0\\
b^2(81\, b^3 + 1000\, \, c^2)=0

\end{array}\]

Si $c=0\longrightarrow 81b^{5}=0\longrightarrow b=0$ car $p\neq 3$.  Si $b=0\longrightarrow 729c^{5}=0$, soit $9^{3}c^{5}=0\longrightarrow c=0$. si $bc\neq 0$ le syst\`{e}me ci-dessus a un d\'{e}terminant nul ~:

\[98.1000-81.729=0\]

Soit $11.3541=0\longrightarrow p=11\ ou \ p=3541$

\begin{itemize}
\item Cas $p=11$ prenons $c=1$ if faut que $b^{3}=-729/98=3 \ mod \ 11\longrightarrow b=-2$, d'o\`{u}~ :

\[P=x^{5}-10x^{3}-3x^{2}\]

On v\'{e}rifie que

\[P=x^{2}(x+1)(x+3)(x-4)\qquad mod\  11\]

On v\'{e}rifie que

\[\begin{array}{l}
P_{2}(-3)=0\  mod\ 11\\
P_{3}(-3)=0\  mod\ 11\\
P_{4}(0)=0\  mod\ 11

\end{array}\]

\item Cas $p=3541$ pour r\'{e}aliser $81\, b^3 + 1000\, \, c^2=0$ il suffit de prendre $b=-10$ et $c=9$~ :

\[P=x^5-10x^3-9x^2\]

On v\'{e}rifie que

\[P=x^{2}(x+1)(x-1567)(x+1566)\qquad mod\  3541\]

On v\'{e}rifie que

\[\begin{array}{l}
P_{2}(1567)=0\  mod\ 3541\\
P_{3}(-1)=0\  mod\ 3541\\
P_{4}(0)=0\  mod\ 3541

\end{array}\]

\end{itemize}

\item 2\`{e}me cas ~:   $b=0$ 

Dans ce cas les resultants  ci-dessus deviennent~:

\[\begin{array}{l}
Res_{x}(P,P_{2})=c^4\cdot \left(100\, a^3 + 729\, c\right)\\
Res_{x}(P,P_{4})=a^2\cdot \left(4\, a^3 + 125\, c\right)

\end{array}\]

Si $P$ v\'{e}rifie les hypoth\`{e}ses de Casas-Alvero on aura alors

\[\begin{array}{l}
c^4\cdot \left(100\, a^3 + 729\, c\right)=0\\
a^2\cdot \left(4\, a^3 + 125\, c\right)=0

\end{array}\]

Si $c=0\longrightarrow 4a^{5}=0\longrightarrow a=0$ car $p\neq 2$.  Si $a=0\longrightarrow 729c^{5}=0$, soit $9^{3}c^{5}=0\longrightarrow c=0$. si $ac\neq 0$ le syst\`{e}me ci-dessus a un d\'{e}terminant nul ~:

\[100.125-4.729=0\]

Soit $9584=2^{4}.599=0\longrightarrow p=599$. Pour r\'{e}aliser $a^2\cdot \left(4\, a^3 + 125\, c\right)=0$ prenons $c=4$ et $a=-5$ d'o\`{u}~ :

\[P=x^{5}+5x^{4}-4x^{2}\]

On v\'{e}rifie que

\[P=x^{2}(x+1)(x+269)(x-265)\qquad mod\  599\]

et

\[\begin{array}{l}
P_{2}(-269)=0\  mod\ 599\\
P_{3}(0)=0\  mod\ 599\\
P_{4}(-1)=0\  mod\ 599

\end{array}\]

\item 3\`{e}me cas ~:   $c=0$ 

Dans ce cas les  resultants ci-dessus deviennent~:

\[\begin{array}{l}
Res_{x}(P,P_{3})=\left(- 9\right)\cdot b^4\cdot \left(2\, a^2 - 9\, b\right)\\
Res_{x}(P,P_{4})=a^3\cdot \left(4\, a^2 - 25\, b\right)

\end{array}\]

Si $P$ v\'{e}rifie les hypoth\`{e}ses de Casas-Alvero on aura alors

\[\begin{array}{l}
\left(- 9\right)\cdot b^4\cdot \left(2\, a^2 - 9\, b\right)=0\\
a^3\cdot \left(4\, a^2 - 25\, b\right)=0

\end{array}\]

Si $b=0\longrightarrow 4a^{5}=0\longrightarrow a=0$ car $p\neq 2$.  Si $a=0\longrightarrow 9^{2}b^{5}=0\longrightarrow b=0$. si $ab\neq 0$ le syst\`{e}me ci-dessus a un d\'{e}terminant nul ~:

\[2.(-25)+4.9=0\]

Soit $-14=-2.7=0\longrightarrow p=7$. Pour r\'{e}aliser $\left(- 9\right)\cdot b^4\cdot \left(2\, a^2 - 9\, b\right)=0$ prenons $b=2$ et $a=3$ d'o\`{u}~ :

\[P=x^{5}-3x^{4}+2x^{3}\]

On v\'{e}rifie que

\[P=x^{3}(x-1)(x-2)\qquad mod\  7\]

et

\[\begin{array}{l}
P_{2}(0)=0\  mod\ 7\\
P_{3}(1)=0\  mod\ 7\\
P_{4}(2)=0\  mod\ 7

\end{array}\]

\end{itemize}

\end{itemize}

\section{Conclusion}
\begin{itemize}
\item  Les resulats ci dessus ne permettent pas de d\'{e}cider  pour les degr\'{e}s : 12,20,24,28,30,35,36,...
\item Les m\'{e}thode ci dessus ne semblent pas s'\'{e}tendre au cas du degr\'{e}s 6, le contre exemple de \cite{Hans}

\[P = X^{6} +3 144 481 702 696 843X^{4} +X^{3} +2 707 944 513 497 181X^{2}\]

\[p=7 390 044 713 023 799\]

 laisse supposer que les mauvais nombres premiers de ce cas sont tr\'{e}s grands et leur nombre est grand

\end{itemize}

\ \\
{\scriptsize
\begin{flushleft}
Polynomials over commutative rings \\
MSC-numbers 2000: 13M10 13P05 13P10 P
\end{flushleft}
\ \\
**************\\
Adresses des auteurs\\
\ \\
**************\\
Prof M.Chellali\\
D\'{e}partement de math\'{e}matiques \\Facult\'{e} des
sciences,  Universit\'{e} Mohammed 1\\Oujda, Maroc.\\
\ \\
**************\\
Prof Alain Salinier\\
D\'{e}partement de Math\'{e}matiques\\
Facult\'{e} des Sciences et Techniques de Limoges\\
123, avenue Albert Thomas\\
87060 LIMOGES Cedex (FRANCE) \\
**************\\
\underline{email}~:
\begin{tabular}{ll}

 & mustapha.chellali@gmail.com\\
 & alain.salinier@unilim.fr
\end{tabular}

\end{document}